\documentclass{amsart}
\usepackage{amsmath, amsfonts, amssymb, amsthm, eufrak}
\usepackage[all]{xy}

\let\k\kappa

\let\s\sigma

\let\f\varphi

\let\Om\Omega

\def\lra{\longrightarrow}

\def\A{\mathcal A}

\def\C{\mathcal C}
\def\E{\mathcal E}
\def\EE{\text{E}}
\def\F{\mathcal F}
\def\G{\mathcal G}
\def\H{\mathcal H}
\def\HH{\text{H}}
\def\hh{\text{h}}

\def\K{\mathcal K}
\def\L{\mathcal L}
\def\LL{\text{L}}
\def\P{\mathbb{P}}
\def\R{\text{R}}

\def\O{\mathcal O}

\def\om{\omega_{{\mathbb P}^n}}
\def\oo{\omega_{S \times {\mathbb P}^n/S}}
\def\ooo{\omega_{R \times {\mathbb P}^n/R}}
\def\M{\text{M}_{{\mathbb P}^n}}
\def\D{{\text{\tiny D}}}
\def\DD{{\text{\emph{\tiny D}}}}

\def\tensor{\otimes}
\def\isom{\simeq}
\def\ext{{\mathcal Ext}}
\def\hom{{\mathcal Hom}}

\def\tor{{\mathcal T\!or}}

\def\Aut{\text{Aut}}

\newcommand{\tilda}{\widetilde}

\def\bdm{\begin{displaymath}}
\def\edm{\end{displaymath}}

\def\ba{\begin{array}}
\def\ea{\end{array}}

\begin{document}

\title[duality for moduli spaces of sheaves supported on projective curves]
{a duality result for moduli spaces of semistable sheaves
supported on projective curves}

\author{mario maican}

\address{Mario Maican \\
Str. Avrig 9-19, Ap. 26 \\
Bucharest, Romania \\
E-mail: m-maican@wiu.edu}

\begin{abstract}
We show that the map sending a sheaf to its dual gives an
isomorphism of the moduli space of semistable sheaves with
fixed multiplicity and Euler characteristic and supported
on projective curves to the moduli space of semistable sheaves 
of dimension one on the projective space with the same
multiplicity but with opposite Euler characteristic.
\end{abstract}

\maketitle

\noindent
Let $k$ be an algebraically closed field of characteristic zero.
Let $\P^n$ be the projective space of dimension $n$ over $k$.
We recall that a coherent algebraic sheaf $\F$ on $\P^n$
has support of dimension one if and only if its Hilbert polynomial
$P_{\F}$ has degree one that is, if we can write $P_{\F}(m)=rm+\chi$.
Here $r=r(\F)$ is a positive integer called the \emph{multiplicity}
of $\F$, while $\chi$ is the Euler characteristic $\chi(\F)=\hh^0(\F)-\hh^1(\F)$.

We fix integers $r \ge 1$ and $\chi$ and we denote by $\M(r,\chi)$
the moduli space of semistable sheaves on $\P^n$ with Hilbert polynomial
$P(m)=rm+\chi$. In this paper we will prove that $\M(r,\chi)$ and $\M(r,-\chi)$
are isomorphic. The isomorphism maps a point represented by $\F$
to the point represented by the dual sheaf $\ext^{n-1}(\F,\om)$.
We will use the concept of semistability due to Gieseker and defined
in terms of the lexicographic order on the coefficients of the reduced
Hilbert polynomial, cf. definition 1.2.4 in \cite{huybrechts}.
We recall that a sheaf $\F$ on $\P^n$ with one-dimensional support
is semistable (stable) if and only if it is pure, meaning that there are no subsheaves
with support of dimension zero, and for every proper subsheaf $\E \subset \F$
we have
\bdm
\frac{\chi(\E)}{r(\E)} \ \le (<) \ \frac{\chi(\F)}{r(\F)}.
\edm

It was first noticed in \cite{fr-diplom}
that $\F \lra \ext^1 (\F,\omega_{\P^2})$
gives a birational map from M$_{\P^2}(r,\chi)$
to M$_{\P^2}(r,-\chi)$ when $r$ and $\chi$ are mutually prime.

We begin by recalling the Beilinson free monad.
Let $\F$ be a sheaf on $\P^n$. A \emph{monad} for $\F$
is a sequence of sheaves
\bdm
0 \lra \C^{-p} \lra \ldots \lra \C^0 \lra \ldots \lra \C^q
\lra 0
\edm
which is exact, except at $\C^0$, where the cohomology is
$\F$. If each $\C^i$ is a direct sum of line bundles
we talk of a \emph{free monad}. In the sequel $\F$ will be
a coherent sheaf on $\P^n$. According to \cite{beilinson},
cf. also \cite{eisenbud}, there is a free monad for $\F$
\bdm
0 \lra \C^{-n} \lra \ldots \lra \C^0 \lra \ldots \lra \C^n
\lra 0
\edm
with
\bdm
\C^i = \bigoplus_{0 \le p \le n} \, \HH^{i+p}(\F \tensor \Om^p(p))
\tensor \O(-p).
\edm

\noindent \\
{\bf Lemma 1:} \emph{Consider a free monad for $\F$}
\bdm
0 \lra \A^{-n} \lra \ldots \lra \A^0 \lra \ldots \lra \A^n
\lra 0
\edm
\emph{with ${\displaystyle \A^i = \oplus_{0 \le p \le n}\,
W_{ip} \tensor \O(-p)}$, where $W_{ip}$ are vector spaces over $k$
of dimension $w_{ip}$.
Then $w_{ip}=\text{\emph{h}}^{i+p}(\F \tensor \Om^p(p))$ for all $i$ and $p$
if and only if each map $W_{i,p} \tensor \O(-p) \lra W_{i+1,p} \tensor \O(-p)$
occuring in the monad is zero.}

\noindent \\
\textsc{Proof:} The two spectral sequences $' \EE$ and $'' \EE$
converging to the hypercohomology of $\A^* \tensor \Om^j(j)$
are given by
\begin{eqnarray*}
' \EE_2^{pq} & = & \HH^p (\P^n, \H^q(\A^* \tensor \Om^j(j))), \\
''\EE_1^{pq} & = & \HH^q (\P^n,\A^p \tensor \Om^j(j)).
\end{eqnarray*}
By hypothesis $\H^q (\A^* \tensor \Om^j(j))$ is
$\F \tensor \Om^j(j)$ for $q=0$ and is zero for $q \neq 0$.
Thus
\bdm
' \EE_2^{pq}=
\begin{cases}
\HH^p(\F \tensor \Om^j(j)), & \text{for  $q=0$}, \\
0                                             & \text{for $q \neq 0$},
\end{cases}
\edm
which shows that $' \EE$ degenerates at level two.
The differentials in the second spectral sequence
$''\EE_1^{p,q} \lra {''\EE}_1^{p+1,q}$ are induced on cohomology
by the maps $\A^p \tensor \Om^j(j) \lra \A^{p+1} \tensor \Om^j(j)$.
According to Bott's formulas on p. 8 in \cite{oss} we have
$\hh^q (\Om^j(j-i)) =0$ for $q \neq j$, $0 \le i \le n$
or for $q=j$, $i \neq j$, $0 \le i \le n$. Thus $''\EE_1^{pq}=0$ if $q \neq j$
and we have the identification
\bdm
''\EE_1^{pj} = W_{pj} \tensor \HH^j(\P^n,\Om^j) \isom W_{pj}.
\edm
The differentials $f_{pj}: {''\EE}_1^{p,j} \lra {''\EE}_1^{p+1,j}$ are identified with
the maps from the monad $W_{p,j} \tensor \O(-j) \lra W_{p+1,j} \tensor \O(-j)$.

Assume now that all maps $f_{pj}$ are zero. Then $''\EE$ degenerates at its first
term and, comparing $'\EE$ with $''\EE$ we get
\bdm
\hh^p(\F \tensor \Om^j(j)) = \text{dim}('\EE_2^{p0}) = \text{dim}('\EE_{\infty}^{p})
= \text{dim}(''\EE_{\infty}^p)=\text{dim}(''\EE_1^{p-j,j})= w_{p-j,j}.
\edm
Conversely, assume that $f_{pj}$ are not all zero. Then $''\EE$ degenerates
at level two and we have
\bdm
\hh^p(\F \tensor \Om^j(j)) = \text{dim}(''\EE_{\infty}^p)=\text{dim}(''\EE_2^{p-j,j})
<\text{dim}(''\EE_1^{p-j,j})= w_{p-j,j}
\edm
for at least one choice of indeces $p$ and $j$. This finishes the proof of the lemma. \\

\noindent \\
{\bf Corollary 2:} \emph{All the maps
$\text{\emph{H}}^i (\F \tensor \Om^p(p)) \tensor \O(-p) \lra
\text{\emph{H}}^{i+1} (\F \tensor \Om^p(p)) \tensor \O(-p)$
occuring in the Beilinson free monad for $\F$ are zero.} \\

In the sequel we will assume that the schematic support
of $\F$ has codimension $c$ in $\P^n$. The \emph{dual sheaf}
$\F^\D$ of $\F$ is defined by $\F^\D = \ext^c (\F,\om)$.
We remark that the hypothesis on the dimension of the support of $\F$
ensures that the extension sheaves $\ext^i(\F,\om)$ vanish for $0 \le i < c$,
see (iii) 7.3 in \cite{hartshorne}.

\noindent \\
{\bf Lemma 3:} \emph{Let $\F$ be a coherent sheaf on $\P^n$
with support of codimension $c \ge 1$. Let}
\bdm
0 \lra \C^{-p} \lra \ldots \lra \C^0 \lra \ldots \lra \C^q
\lra 0
\edm
\emph{be a free monad for $\F$. Assume that $\ext^i (\F,\om)=0$
for $i > c$. We consider the dual bundles
$\C^i_\DD = \hom (\C^{-i-c},\om)$. Then the dual sequence}
\bdm
0 \lra \C_\D^{-q-c} \lra \ldots \lra \C_\D^0 \lra \ldots
\lra \C_\D^{p-c} \lra 0
\edm
\emph{is a free monad for $\F^\DD$.}

\noindent \\
\textsc{Proof:} We compare the spectral sequences $'\EE$ and $''\EE$
converging to the hyper-derived functor $\ext^{p+q}(\C^{\bullet},\om)$.
They are given by
\begin{align*}
'\EE_1^{pq} & = \ext^q (\C^{-p},\om), \\
''\EE_2^{pq} & = \ext^p (\H^{-q}(\C^*),\om).
\end{align*}
Since $\C^i$ are projective sheaves, we have $'\EE_1^{pq} = 0$ for $q \neq 0$.
This shows that $'\EE$ degenerates at level two.
The same is true of $''\EE$ because $''\EE_2^{pq}=0$ for $q \neq 0$.
By virtue of the hypothesis and of the remark preceding the lemma,
we have
\bdm
''\EE_2^{p0} = \begin{cases}
\F^\D & \text{for $p=c$}, \\
0 & \text{for $p \neq c$}.
\end{cases}
\edm
In conclusion,  the zero row of $'\EE_1$ provides the desired monad for $\F^\D$. \\

For any sheaf $\F$ on $\P^n$ there is a natural homomorphism
$\F \lra \F^{\D\D}$ which is injective if and only if $\F$ is pure.
We recall that a sheaf with support of dimension $d$ is called \emph{pure}
if it does not have a
nonzero subsheaf with support of dimension smaller than $d$.
We say that $\F$ is \emph{reflexive} if the map $\F \lra \F^{\D\D}$
is an isomorphism.
According to 1.1.10 in \cite{huybrechts}
the hypotheses of lemma 3 are satisfied for pure sheaves
of dimension one and for reflexive sheaves of dimension two. 

\noindent \\
{\bf Remark 4:} Let $\F$ be a sheaf on $\P^n$ with support
of dimension one. We assume that $\F$ is pure and we notice
that this is equivalent to saying that $\F$ has no
zero-dimensional torsion. According to lemma 3.1(i) and proposition 3.3(iv)
from \cite{sga2},
this is further equivalent to saying that at every closed
point $x$ in the support of $\F$ we have depth$(\F_x) \ge 1$.
From this we see that $\F$ satisfies Serre's
condition S$_{2,n-1}$:
\bdm
\text{depth}(\F_x) \ge \text{min} \{ 2,\ \text{dim}(\O_{\P^n,x})
-n+1 \} \text{ for all } x \in \text{Supp}(\F).
\edm
From 1.1.10 in \cite{huybrechts} we conclude that $\F$ is
reflexive.

\noindent \\
{\bf Proposition 5:} \emph{Let $\F$ be a coherent sheaf on
$\P^n$ with support of codimension $n-d=$ $c \ge 1$. Assume that
$\ext^i (\F,\om)=0$ for $i > c$. Then for all $i$, $j$ we have}
\bdm
\hh^i (\F \tensor \Om^j(j))= \hh^{d-i} (\F^\D \tensor \Om^{n-j}(n-j+1)).
\edm

\noindent \\
\textsc{Proof:} We apply lemma 3 to the Beilinson free monad
for $\F$. We get a monad for $\F^\D$ with terms
\bdm
\C^i_\D = \hom (\C^{-i-c},\om) =
\bigoplus_{0 \le p \le n} \, \HH^{-i-c+p}(\F \tensor \Om^p(p))
\tensor \O(p-n-1).
\edm
We tensor the above monad with $\O(1)$ to get a monad for $\F^\D(1)$.
In view of corollary 2, this monad satisfies the necessary and sufficient
condition from lemma 1. We conclude that for all $i$, $p$
\bdm
\hh^{i+p}(\F^\D \tensor \Om^p (p+1)) =
\hh^{d-i-p} (\F \tensor \Om^{n-p}(n-p)).
\edm
This proves the proposition.

\noindent \\
{\bf Corollary 6:} \emph{Let $\F$ be as in the previous
proposition. Then for all $i$, $j$ we have}
\bdm
\hh^i (\F) = \hh^{d-i}(\F^\D).
\edm

\noindent \\
{\bf Remark 7:} The above corollary also follows from Serre duality and the
degeneration of the local-to-global spectral sequence
$\EE_2^{pq}= \HH^p (\P^n, \ext^q (\F,\om))$,
which converges to $\text{Ext}^{p+q}(\F,\om)$. \\

Next we would like to relate the Hilbert polynomials $P_{\F}$
and $P_{\F^\D}$ of $\F$, respectively $\F^\D$. We recall that
$P_{\F}$ has degree equal to the dimension of Supp$(\F)$.
As the support of $\F^\D$ is included in the support of $\F$,
we have deg$(P_{\F^\D}) \le $ deg$(P_{\F})$.

\noindent \\
{\bf Corollary 8:} \emph{Let $\F$ be as in the previous
proposition. Then for all $m$ we have} 
\bdm
P_{\F^\D}(m)=(-1)^{d} P_{\F}(-m).
\edm
\emph{In particular, if $\F$ is pure of dimension one with
Hilbert polynomial $P_{\F}(m)=rm+\chi$, then the dual sheaf has Hilbert polynomial
$P_{\F^\DD}(m)= rm-\chi$.}

\noindent \\
{\bf Lemma 9:} \emph{Let $\F$ be a coherent sheaf on $\P^n$,
$n \ge 2$, pure of dimension one. Then $\F$ is semistable (stable) if and
only if $\F^\DD$ is semistable (stable).}

\noindent \\
\textsc{Proof:} According to 1.2.6 in \cite{huybrechts}, the sheaf $\F$ is semistable
(stable) if and only if for all pure one-dimensional quotients $\G$ of $\F$ we have
\bdm
\frac{\chi(\G)}{r(\G)} \ \ge (>) \ \frac{\chi(\F)}{r(\F)}.
\edm
Take a pure one-dimensional destabilizing quotient $\G = \F/\K$.
Since $\K$ has support of dimension one, we have $\ext^{n-2}(\K,\om)=0$, hence
$\G^\D$ is a subsheaf of $\F^\D$. From corollary 8 we get
\bdm
\frac{\chi(\G^\D)}{r(\G^\D)}=-\frac{\chi(\G)}{r(\G)}>-\frac{\chi(\F)}{r(\F)}=\frac{\chi(\F^\D)}{r(\F^\D)}.
\edm
Thus $\G^\D$ is a destabilizing subsheaf of $\F^\D$. This proves sufficiency.
Necessity follows from the
fact that $\F$ is reflexive, cf. remark 4. \\

Two semistable sheaves $\F$ and $\G$ on $\P^n$ with Hilbert polynomial $P$
give the same point in $\M(P)$ if and only if they are \emph{S-equivalent},
meaning that there are filtrations by subsheaves
\bdm
0 = \F_0 \subset \F_1 \subset \ldots \subset \F_{\k-1} \subset \F_{\k} = \F,
\edm
\bdm
0 = \G_0 \subset \G_1 \subset \ldots \subset \G_{\k-1} \subset \G_{\k} = \G,
\edm
such that all quotients $\F_i/\F_{i-1}$ and $\G_i/\G_{i-1}$ are stable and there is
a permutation $\s$ of the set of indeces $\{ 1, \ldots, \k \}$ such that for all $i$
\bdm
\F_i/\F_{i-1} \isom \G_{\s(i)}/\G_{\s(i)-1}.
\edm
For stable sheaves S-equivalence means isomorphism.
In other words, if the above conditions are satisfied, then $\F$ is stable if and only if
$\G$ is stable and then $\k=1$ and $\F \isom \G$. A \emph{Jordan-H\"older filtration}
of a semistable sheaf is a filtration as above by subsheaves such that all quotients
are stable. \\

As in the proof of 9.3 from \cite{maican}, a Jordan-H\"older filtration
\bdm
0 = \F_0 \subset \F_1 \subset \ldots \subset \F_{\k-1} \subset \F_{\k} = \F
\edm
for a semistable sheaf of dimension one on $\P^n$ gives a Jordan-H\"older filtration
\bdm
0 = (\F/\F_{\k})^\D \subset (\F/\F_{{\k}-1})^\D \subset \ldots \subset
(\F/\F_1)^\D \subset (\F/\F_0)^\D = \F^\D
\edm
for the dual sheaf with quotients $(\F_{i}/\F_{i-1})^\D$.
The latter are stable by virtue of lemma 9. We arrive at the following lemma:

\noindent \\
{\bf Lemma 10:} \emph{Let $\F$ and $\G$ be semistable sheaves on $\P^n$
with the same Hilbert polynomial $P(m)=rm+\chi$.
If $\F$ and $\G$ are S-equivalent, then so are their duals.}

\noindent \\
{\bf Lemma 11:} \emph{Let $S$ be an algebraic scheme over $k$
and let $\F$ be an $S$-flat coherent sheaf on $S \times \P^n$.
Assume that for every $s$ in $S$ the restriction $\F_s = \F_{|\{ s \} \times \P^n}$
has support of codimension $c$. Then for all $i < c$ we have}
\bdm
\ext^i (\F,\oo)=0.
\edm

\noindent \\
\textsc{Proof:} Without loss of generality we may assume that $S$
is affine. As in the proof of (iii) 7.3 from \cite{hartshorne},
we reduce to showing the vanishing of Ext$^i(\F,\oo(q))$
for large $q$. By 11.2(f) on p. 213 of \cite{residues}
we have the duality
\bdm
\text{Ext}_{\O_{S \times \P^n}}^i (\F,\oo(q)) \isom \text{Hom}_{\O_S}(\R^{n-i}f_* (\F(-q)),\O_S)
\edm
where $f: S \times \P^n \lra S$ is the projection onto the first component.
As $n-i$ exceeds the dimension of each restriction $\F(-q)_s$, we have
$\HH^{n-i}(\F(-q)_s)=0$ for all $s$ in $S$.
From exercise (iii) 11.8 in \cite{hartshorne} we deduce that $\R^{n-i}f_* (\F(-q))=0$,
so the above extension group vanishes. \\

In the sequel we will need the following relative version of the Hilbert syzygy theorem.
Let $S$ be an algebraic scheme over $k$ and $\F$ a coherent $S$-flat sheaf
on $S \times \P^n$. Then there exists a locally free resolution $\E^{^{\bullet}} \lra \F$
of length at most $n$. Such a resolution can be constructed using the
relative Beilinson spectral sequence with first term
\bdm
\EE_1^{pq}= \R^q f_* (\F(m) \tensor g^* \Om^{-p}(-p)) \boxtimes \O_{\P^n}(p),
\edm
which converges to $\F(m)$ in degree zero and to 0 in degree different from zero
(see 4.1.11 in \cite{oss}).
Here $f : S \times \P^n \lra S$ and $g : S \times \P^n \lra \P^n$ are the canonical
projections. If $m$ is large enough, then $\EE_1^{pq}=0$ for $q \neq 0$ and
the push-forward sheaves $f_* (\F(m) \tensor g^* \Om^{-p}(-p))$ are locally free.
Thus the spectral sequence degenerates at level two
and the zero row of $\EE_1$ provides a locally free resolution of $\F(m)$.
We can cover $S$ with open affine subsets $U$, such that the restriction
of each term of the resolution to $U \times \P^n$ is a direct sum of line bundles
of the form $g^* \O_{\P^n}(-l)$.

\noindent \\
{\bf Lemma 12:} \emph{Let $S$ be an algebraic scheme over $k$ and
let $\F$ be an $S$-flat coherent sheaf on $S \times \P^n$.
Assume that for every $s$ in $S$
the restriction $\F_s = \F_{|\{ s \} \times \P^n}$ has support
of codimension $c$ and that $\ext^i (\F_s,\om)=0$ for $i>c$.
Then}
\bdm
\ext^i(\F,\oo) = 0 \quad \text{\emph{for}} \quad i \neq c.
\edm
\emph{Moreover, the sheaf $\F^\DD=\ext^c(\F,\oo)$ is flat over $S$ and
$(\F_s)^\DD \isom (\F^\DD)_s$ for all $s$ in $S$.}

\noindent \\
\textsc{Proof:} Consider a finite locally free resolution $\E^{^{\bullet}} \lra \F$.
The extension sheaf $\ext^{p}(\F, \oo)$ is the $p$-th cohomology of
the complex $\hom(\E^{^{\bullet}},\oo)$. We fix a point $s$ in $S$ and denote
$i : \{ s \} \times \P^n \lra S \times \P^n$ the canonical inclusion.
The functor $\hom(\_\!{\hspace{0.3mm}}\_, \oo)$ sends projective sheaves to $i^*$-acyclic sheaves,
hence there is the Grothendieck spectral sequence (see viii 9.3 in \cite{hilton-stammbach})
\bdm
\EE_2^{pq} = \LL^p i^* \, \R^q \hom(\_\!\hspace{0.3mm}\_,\oo)(\F)
= \tor_{-p}^{\O_{S \times \P^n}} (\ext^q_{\O_{S \times \P^n}} (\F, \oo), \O_{\{s \} \times \P^n}),
\edm
which converges to the hyper-derived functors
${\mathbb L}^{p+q}i^*(\hom(\E^{^{\bullet}},\oo))$.
We recall that the above is the second term of the spectral sequence
associated to the bigraded complex obtained by taking a fourth-quadrant
Cartan-Eilenberg resolution of $\hom(\E^{^{\bullet}},\oo)$ and then applying $i^*$.
The other spectral sequence associated to the bigraded complex is given by
\bdm
'\EE_1^{pq}= \LL^q i^*(\hom(\E^{-p},\oo)).
\edm
$'\EE$ degenerates at level two because $'\EE_1^{pq}=0$ for $q \neq 0$. Thus
\bdm
{\mathbb L}^p i^* (\hom(\E^{^{\bullet}},\oo)) \isom {'\EE}_{\infty}^{p0} = {'\EE}_2^{p0}
\edm
is the cohomology at position $i^* \hom(\E^{-p},\oo)$ of the complex
\bdm
i^* \hom(\E^{^{\bullet}},\oo) \isom \hom(i^*\E^{^{\bullet}},i^*\oo) \isom \hom(\E_s^{^{\bullet}},\om).
\edm
As every term of the resolution $\E^{^{\bullet}}\lra \F$ is $S$-flat,
restricting to $\{ s \} \times \P^n$ we obtain a locally free resolution
$\E^{^{\bullet}}_s \lra \F_s$. Thus, for all $j$,
\bdm
{\mathbb L}^j i^*(\hom(\E^{^{\bullet}},\oo)) \isom \ext^j(\F_s,\om).
\edm
By hypothesis the above sheaves are zero for $j \neq c$,
forcing $\EE_{\infty}^{pq}=0$ for $p+q \neq c$.
Assume $c<n$. From the vanishing of
\bdm
\EE_{\infty}^{0n} = \EE_2^{0n} = \ext^n(\F,\oo)_s
\edm
for an arbitrary point $s$ in $S$, we deduce that $\ext^n(\F,\oo)=0$.
Thus $\EE_2^{pn}=0$ for all $p$, forcing $\EE_2^{0,n-1}= \EE_{\infty}^{0,n-1}$.
If $c < n-1$ we can repeat the argument and conclude that
$\ext^{n-1}(\F,\oo)=0$. By induction we obtain that the sheaves
$\ext^q(\F,\oo)$ vanish for $q>c$. In view of lemma 11 they also vanish for $q < c$.
Thus $\EE_2^{pq}=0$ for $q \neq c$ and $\EE$ degenerates at level two:
\bdm
\EE_2^{pc} = \EE_{\infty}^{pc} \isom {\mathbb L}^{p+c}i^*(\hom(\E^{^{\bullet}},\om))
\isom \ext^{p+c}(\F_s,\om).
\edm
For $p=0$ we arrive at the isomorphism $(\F^\D)_s \isom (\F_s)^\D$.
For $p=-1$ we get
\bdm
\tor_1^{\O_{S \times \P^n}} (\F^\D, \O_{\{ s \} \times \P^n}) = 0.
\edm
From this we deduce that $\F^\D$ is $S$-flat. Indeed, localizing at a point $x$ lying over $s$,
we obtain
\bdm
\text{Tor}_1^{\O_x}((\F^\D)_x, \O_x \tensor_{\O_s} k(s)) =0.
\edm
Here $k(s)$ denotes the residue field of $s$. We now use the local criterion of flatness
4.1.3 in \cite{lepotier} to conclude that the stalk $(\F^\D)_x$ is flat over $\O_s$. \\

Before we proceed any further we need to review the construction of the moduli
space $\M (P)$ of semistable sheaves on $\P^n$ with fixed
Hilbert polynomial $P$. There exists
an integer $m \gg 0$ such that for every semistable sheaf $\F$ on $\P^n$
with Hilbert polynomial $P$ the twisted sheaf $\F(m)$ is generated by global sections
and the higher cohomology groups $\HH^i(\F(m))$, $i \ge 1$, vanish.
Thus $\hh^0(\F(m))=P(m)$ and $\F$ occurs as a quotient
$\rho : V \tensor \O(-m) \twoheadrightarrow \F$, where $V$ is a fixed vector space over $k$
of dimension $P(m)$.
We consider the quotient scheme $Q = Quot_{\P^n}(V \tensor \O_{\P^n}(-m),P)$
and the universal quotient sheaf $\tilda{\F}$ on $Q \times \P^n$.
Let $R \subset Q$ be the open subset of equivalence classes of quotients
$[ \rho : V \tensor \O_{\P^n}(-m) \twoheadrightarrow \F]$
for which $\F$ is semistable and the map on global sections
$\HH^0(\rho(m)): V \lra \HH^0(\F(m))$ is an isomophism.
The reductive group SL$(V)$ acts on $Q$ via its action on the first
component of $V \tensor \O_{\P^n}(-m)$. Clearly $R$ is invariant
under this action. If $m$ is chosen large enough, then
$\M (P)$ is a good quotient of $R$ by SL$(V)$. We denote by
$\pi : R \lra \M (P)$ be the good quotient map.

\noindent \\
{\bf Theorem 13:} \emph{For any integers $n \ge 2$, $r \ge 1$ and
$\chi$ there is an isomorphism}
\bdm
\M (r,\chi)\lra \M (r,-\chi)
\edm
\emph{which maps the S-equivalence class of a sheaf $\F$
to the S-equivalence class of $\F^\DD$.}

\noindent \\
\textsc{Proof:} According to remark 4, corollary 8, lemma 9 and lemma 10
the above map, call it $\eta$, is well-defined and bijective. It remains to
show that $\eta$ is a morphism. By symmetry it will follow
that $\eta^{-1}$ is also a morphism.

We put $P(m)=rm+\chi$, so $P^\D (m)=rm-\chi$. We adopt the
notations preceeding the theorem. The sheaf $\tilda{\F}^\D$
on $R \times \P^n$ defined by
$\tilda{\F}^\D = \ext^{n-1} (\tilda{\F},\ooo)$
is coherent and, by the previous lemma, $R$-flat.
Moreover, for each $s$ in $R$,
its restriction $\tilda{\F}^\D_s$ to $\{ s \} \times \P^n$
is isomorphic to $(\tilda{\F}_s)^\D$. According to corollary 8,
the latter has Hilbert polynomial $P^\D$.
The moduli space property of $\M (P^\D)$ gives
a morphism $\pi^\D :R \lra \M (P^\D)$ which maps $s$ to the
S-equivalence class of $\tilda{\F}^\D_s$.

From lemma 10 we see that $\pi^\D$ is constant on the fibers of
$\pi$. The good quotient property of $\pi$ shows that
$\pi^\D$ factors through a morphism $\M (P) \lra \M (P^\D)$.
This morphism is $\eta$, which finishes the proof of the
theorem. \\

\noindent
{\bf Remark 14:} If $P$ has degree at least 2, then the above proof yields
a morphism
\bdm
U \lra \M (P^\D), \qquad [\F] \lra [\F^\D],
\edm
defined on the open set $U$ in $\M (P)$ of isomorphism classes
of stable sheaves $\F$ with $\F^\D$ stable and satisfying
$\ext^i (\F,\om)=0$ for $i>c$. For sheaves supported on surfaces the last
condition is equivalent to saying that $\F$ is reflexive.
We notice that $U$ is nonempty for all $P$ of degree 2
for which there is a smooth surface $X$ in $\P^n$ and a line
bundle ${\mathcal L}$ on $X$ with Hilbert polynomial $P$.
This is so because any locally free $\O_X$-module of rank one,
with $X$ a reduced subscheme of $\P^n$,
is stable as a sheaf on $\P^n$.
Moreover, $\L$ is reflexive as a sheaf on $\P^n$, which can be seen from
the isomorphism $(i_* {\mathcal L})^\D \isom i_* ({\mathcal L}^\D)$,
where $i : X \lra \P^n$ is the embedding map. \\

In the remaining part of this paper we will apply the duality result to the
following situation. Assume that $\M(r,\chi)$ is the quotient of a certain
parameter space by an algebraic group. This could be the case, for instance,
if all sheaves giving a point in $\M(r,\chi)$ are quasi-isomorphic to
monads of a certain kind. Constructions of moduli spaces of sheaves as quotients
of parameter spaces are ubiquitous in the literature, mostly involving actions
of reductive groups, and, more recently, as in \cite{drezet-1991},
\cite{drezet-trautmann} or \cite{fr-trautmann} for actions of nonreductive groups.
Under the above hypothesis on $\M(r,\chi)$, we will show that $\M(r,-\chi)$
is the quotient of the ``dual parameter space'' modulo the ``dual group.''
This will generalize our partial results from section 9 of \cite{maican}, where
we dealt with locally closed subvarieties inside $\text{M}_{\P^2}(r,\chi)$
for certain choices of $r$ and $\chi$.

We fix sheaves $\C^i$ on $\P^n$, $-p \le i \le q$,
that are finite direct sums of line bundles.
We fix a polynomial $P(m)=rm+\chi$.
We assume that each semistable sheaf
on $\P^n$ with Hilbert polynomial $P(m)=rm+\chi$ occurs
as the cohomology of a monad
\bdm
0 \lra \C^{-p} \stackrel{\f_{-p}\,}{\lra} \ldots
\stackrel{\f_{-1}\,}{\lra} \C^0 \stackrel{\f_0}{\lra} \ldots
\stackrel{\f_{q-1}}{\lra} \C^q \lra 0.
\edm
For a monad of this form we write $\f=(\f_{-p},\ldots,\f_{q-1})$
and we denote by $\F_{\f}$ its cohomology.
Let $W^{ss}$ be the set of all those $\f$ for which $\F_{\f}$
is semistable. We assume that $W^{ss}$ is a constructible subset
inside the finite dimensional vector space
\bdm
W= \text{Hom}(\C^{-p},\C^{-p+1}) \times \ldots \times
\text{Hom}(\C^{q-1},\C^q).
\edm
We equip $W^{ss}$ with the induced structure of a reduced variety.
The algebraic group
\bdm
G = \text{Aut}(\C^{-p}) \times \ldots \times \text{Aut}(\C^q)
\edm
acts on $W$ in an obvious manner: given $g=(g_{-p},\ldots,g_q)$
and a tuple $\f$ as above we put
\bdm
g \, \f =
(g_{-p+1}\f_{-p} \, g_{-p}^{-1}, \ldots, g_{q} \f_{q-1} g_{q-1}^{-1}).
\edm
We remark that, if at least one $\C^i$
has a direct summand of the form $\O_{\P^n}(a) \oplus \O_{\P^n}(b)$,
$a < b$, then $G$ is nonreductive, i.e. it has a nontrivial unipotent radical,
cf. 19.5 in \cite{humphreys}. Indeed, the subgroup of $\Aut(\O_{\P^n}(a) \oplus
\O_{\P^n}(b))$ given by matrices of the form
\bdm
\left[
\begin{array}{cc}
1 & 0 \\
\star & 1
\end{array}
\right]
\edm
is a normal connected unipotent subgroup. This shows that the unipotent
radical of $\Aut(\O_{\P^n}(a) \oplus \O_{\P^n}(b))$ is nontrivial.
A similar argument applies in general to $G$.

Clearly $\F_{\f}$ and $\F_{g \f}$ are isomorphic,
so one is semistable if and only if the other is.
This shows that $W^{ss}$ is $G$-invariant.

Let $p : W^{ss} \times \P^n \lra \P^n$ be the projection onto the
second component. We put $\tilda{\C}^i = p^* \C^i$.
On $W^{ss} \times \P^n$ we have the universal monad
\bdm
\tilda{\C}^{-p} \stackrel{\Phi_{-p}\,}{\lra} \ldots
\stackrel{\Phi_{-1}\,}{\lra} \tilda{\C}^0
\stackrel{\Phi_0}{\lra} \ldots
\stackrel{\Phi_{q-1}}{\lra} \tilda{\C}^q
\edm
with ${\Phi_i}_{|\{ \f \} \times \P^n}=\f_i$.
Let $\tilda{\F}$ denote the middle cohomology of this sequence.
Using arguments as in the proof of lemma 12 one can see that
its restriction $\tilda{\F}_{\f}$ to any fiber $\{ \f \} \times \P^n$
is isomorphic to $\F_{\f}$. In particular, all $\tilda{\F}_{\f}$
have the same Hilbert polynomial $P$. As $W^{ss}$ is reduced
we deduce that $\tilda{\F}$ is $W^{ss}$-flat.
The moduli space property of $\M (P)$ yields a morphism
\bdm
\pi : W^{ss} \lra \M (P), \qquad \pi (\f)=[\F_{\f}].
\edm
Clearly $\pi$ is $G$-equivariant. We will be interested in the
case in which $\pi$ is a good or a geometric quotient. \\

We now apply the duality result to the above set-up.
Under the hypothesis that each sheaf
giving a point in $\M(P)$ is the cohomology of a monad $\C^*$ and by virtue
of lemma 3, corollary 8 and lemma 9, we see that each semistable sheaf on $\P^n$
with Hilbert polynomial $P^\D (m) = rm-\chi$ is the cohomology of a monad
\bdm
0 \lra \C^{-q-c}_\D \stackrel{\f^\D_{-q-c}}{\lra} \ldots
\stackrel{\f^\D_{-1}}{\lra} \C^0_\D \stackrel{\f^\D_0}{\lra} \ldots
\stackrel{\f^\D_{p-c-1}}{\lra} \C_\D^{p-c} \lra 0
\edm
with
\bdm
\C^i_\D = \hom (\C^{-i-c},\om), \qquad \f_i^\D = \hom (\f_{-i-c-1},\om).
\edm
We let $W^{ss}_\D$ be the set of all those tuples
$\f^\D=(\f^\D_{-q-c}, \ldots, \f^\D_{p-c-1})$ for which
$\F_{\f^\D}$ is semistable. We view $W^{ss}_\D$ as a subvariety inside
the affine variety
\bdm
W_\D =
\text{Hom}(\C_\D^{-q-c},\C_\D^{-q-c+1})
\times \ldots \times \text{Hom}(\C_\D^{p-c-1},\C_\D^{p-c})
\edm
equipped with the analogously defined action of the group of automorphisms
\bdm
G_\D = \Aut(\C_\D^{-q-c}) \times \ldots \times \Aut(\C_\D^{p-c}).
\edm
As before, there is a $G_\D$-equivariant map
\bdm
\pi^\D : W^{ss}_\D \lra \M (P^\D), \qquad \pi^\D (\f^\D)=[\F_{\f^\D}].
\edm

\noindent \\
{\bf Corollary 15:} \emph{$\pi$ is a good (geometric) quotient
map if and only if $\pi^\D$ is a good (geometric) quotient map.}

\noindent \\
\textsc{Proof:} First we notice that $G \isom G_\D$ and that the map
$\f \lra \f^\D$ acts by transposition, so it gives
an equivariant isomorphism $W^{ss} \lra W^{ss}_\D$.
According to lemma 3 this isomorphism fits into a commutative diagram
\bdm
\xymatrix
{
W^{ss} \ar[r] \ar[d]_{\pi} & W^{ss}_\D \ar[d]^{\pi^\D} \\
\M(P) \ar[r] & \M (P^\D)
}.
\edm
By virtue of the previous theorem the bottom map $[\F] \lra [\F^\D]$
is an isomorphism, which proves the statement. \\

We finish with an example. Fix a vector space $V$ over $k$ of dimension 4. 
According to \cite{fr-trautmann} a sheaf $\F$ on $\P^3=\P(V)$
with Hilbert polynomial $3m+1$ is semistable if and only if it has a resolution
\bdm
0 \lra 2 \O(-3) \stackrel{\psi}{\lra} \O(-1) \oplus 3\O(-2)
\stackrel{\f}{\lra} \O \oplus \O(-1) \lra \F \lra 0
\edm
with $\f$ not equivalent to a matrix of the form
\bdm
\left[
\ba{cccc}
\star & \star & \star & \star \\
0 & 0 & \star & \star
\ea
\right] .
\edm
Moreover, M$_{\P^3}(3,1)$ is a geometric quotient of the space
of parameters $(\psi,\f)$ modulo the action of the group
of automorphisms. The concrete description of the space of parameters $W^{ss}$
is given at 5.2 in \cite{fr-trautmann}. It consists of pairs $(\psi,\f)$ for which the sequence
of global sections
\bdm
0 \lra k^2 \stackrel{\HH^0(\psi(3))}{\lra} S^2 V^* \oplus (k^3 \tensor V^*)
\stackrel{\HH^0(\f(3))}{\lra} S^3 V^* \oplus S^2 V^*
\edm
is exact and either $\f_{21} \neq 0$ or the entries $\{ \f_{22}, \f_{23}, \f_{24} \}$
are linearly independent as elements of $V^*$.

The above considerations tell us that each
sheaf $\F$ on $\P^3$ with Hilbert polynomial $3m-1$ is semistable if and only if
there is $(\psi,\f)$ in $W^{ss}$ and a resolution
\bdm
0 \lra \O(-4) \oplus \O(-3) \stackrel{\f^\D}{\lra}
\O(-3) \oplus 3\O(-2) \stackrel{\psi^\D}{\lra} 2\O(-1) \lra \F
\lra 0.
\edm
From the above corollary we deduce that $\text{M}_{\P^3}(3,-1)$ is a
geometric quotient of the space of parameters $W^{ss}_\D$ of transposed
matrices $(\f^\D,\psi^\D)$ modulo the group of automorphisms. \\

\begin{center}
\textsc{Acknowledgements}
\end{center}
\vspace{1mm}

\noindent
The author thanks the referee for providing several corrections and improvements
and for numerous helpful comments. The referee suggested the proof of lemma 12,
shortened our original proof of lemma 3 and pointed out that
the two conditions from lemma 1 are equivalent (originally we stated only sufficiency).

%%%%%%%%%%%%%%%%%%%%%%%%%%%%%%%%%%%%%%%%%%%%%%%%%%%%%%%%%%%%%%%%%%

\end{document}